\documentclass{amsart}

\title[Reverse Littlewood--Paley, restriction and Kakeya]{A remark on reverse Littlewood--Paley, restriction and Kakeya}
\author{Anthony Carbery}
\thanks{}
\address{Anthony Carbery, 
School of Mathematics and Maxwell Institute for Mathematical Sciences, 
University of Edinburgh,
JCMB, 
King's Buildings, 
Mayfield Road, 
Edinburgh, EH9 3JZ, 
Scotland.} 
\email{A.Carbery@ed.ac.uk}

\date{June 2014}
\usepackage{amssymb}
\usepackage{amsmath}
\usepackage{url}

\newtheorem{proposition}{Proposition}

\begin{document}

\begin{abstract}
We show that a certain conjectured optimal reverse Littlewood--Paley inequality would, if true, imply
sharp results for the Kakeya maximal function, the Bochner--Riesz means and the Fourier 
restriction operator.
\end{abstract}

\maketitle
\section{Introduction}
\noindent
Let $\delta >0$ be a small parameter and let $\Phi: \mathbb{R} \to \mathbb{R}$ be a smooth function of 
compact support in $[-1, 1]$ satisfying $|\Phi^{(k)}(t)| \leq C_k$ for all $k \in \mathbb{N}$. Define
the Fourier multiplier operator $S^\delta$ on $\mathbb{R}^n$ by 
$$ \widehat{(S^\delta f)}(\xi) = \Phi\left(\frac{|\xi| -1}{\delta}\right) \widehat{f}(\xi).$$
We decompose the $\delta$-neighbourhood $\{ \xi \, : \,  \big| |\xi| - 1 \big| \leq \delta \}$ of the unit sphere 
$\mathbb{S}^{n-1}$ into coin-shaped pieces $E_\alpha$ of tangential dimensions 
$\delta^{1/2} \times \dots \times \delta^{1/2}$ and radial dimension $\delta$, 
and correspondingly the operator 
$S^\delta = \sum_\alpha S_\alpha$ where $S_\alpha$ is a Fourier multiplier 
operator with smooth multiplier $\phi_\alpha$ supported in and adapted to $E_\alpha$. We consider
the reverse Littlewood--Paley inequality
\begin{equation}\label{RLP}
\|S^\delta f \|_{L^{\frac{2n}{n-1}}(\mathbb{R}^n)} \leq 
C_n \left\| \left(\sum_\alpha |S_\alpha f|^2 \right)^{1/2} \right\|_{L^{\frac{2n}{n-1}}(\mathbb{R}^n)}
\end{equation}
where $C_n$ is supposed to be independent of $\delta$ and $f$. When $n=2$ (and thus $2n/(n-1) = 4$), 
Fefferman in \cite{F} proved that this inequality is true since we can multiply out the $L^4$ norm 
and observe that the algebraic differences $ E_\alpha -  E_{\alpha^\prime}$ are essentially disjoint as 
$\alpha \neq \alpha^\prime$ vary. In higher dimensions it remains open. 

\medskip
\noindent
More generally one might ask whether for $q \geq 2$ we have
$$\|S^\delta f \|_{L^{q}(\mathbb{R}^n)} \leq 
C_n \left\| \left(\sum_\alpha |S_\alpha f|^r \right)^{1/r} \right\|_{L^{q}(\mathbb{R}^n)},$$
where for $2 \leq q \leq 2n/(n-1)$ we take $r=2$ and when $q \geq 2n/(n-1)$ we take
$r^\prime = q(n-1)/n$ (so that when $q = \infty$ we have $r=1$). The issue of such inequalities 
with the order the mixed norms reversed, first proposed by Bonami and Garrig\'os, has recently been 
studied in \cite{BD} where they are termed $l^r$-decoupling inequalities.
  
\medskip
\noindent
There is a maximal function relevant to the study of $S^\delta$, the so-called Nikodym maximal function,
for which Fefferman also proved in \cite{F} optimal $L^2$ bounds in two dimensions (in slightly disguised form). 
By what are very familiar arguments (see for example \cite{CoD}, \cite{CoTIO}) these two ingredients can be combined to 
prove the optimal two-dimensional result
$$ \|S^\delta f \|_4 \leq C \left( \log \frac{1}{\delta}\right)^{1/4} \|f \|_4$$
for the operators $S^\delta$.

\medskip
\noindent
In this note we show that (a variant of) \eqref{RLP} actually implies the correct $L^n$ bound for the maximal function, and thus 
\eqref{RLP} together with its variant gives the optimal Bochner--Riesz multiplier result
\begin{equation}\label{BR}
\|S^\delta f \|_{L^{\frac{2n}{n-1}}(\mathbb{R}^n)} \leq
C_n \left( \log \frac{1}{\delta}\right)^{\frac{n-1}{2n}} \|f \|_{L^{\frac{2n}{n-1}}(\mathbb{R}^n)}
\end{equation}
in all dimensions upon combining \eqref{RLP} with the maximal function estimate in the familiar way.

\medskip
\noindent
For ease of exposition we choose to work in the alternative and essentially equivalent setting 
of the extension problem for the Fourier transform. Since the kernels of $S^\delta$ and $S_\alpha$ in 
\eqref{RLP} are essentially localised at scale $1/\delta$, \eqref{RLP} is equivalent to the corresponding 
local inequality
\begin{equation}\label{RLPloc}
\|S^\delta f \|_{L^{\frac{2n}{n-1}}(B(0,\delta^{-1}))} \leq 
C_n \left\| \left(\sum_\alpha |S_\alpha f|^2 \right)^{1/2} \right\|_{L^{\frac{2n}{n-1}}(B(0,\delta^{-1}))}
\end{equation}
and thus to\footnote{Strictly speaking \eqref{RLP} is equivalent to the variant of \eqref{mainlocal} 
where we {\em average} over spheres of radius $r \in [1- \delta, 1+ \delta]$ on both sides, ensuring
that we may assume that the integrand on the right hand side is indeed essentially supported in 
$B(0,\delta^{-1})$. Without this interpretation it is not clear what meaning \eqref{mainlocal} may have 
in general. However, as we shall see below, we shall actually use \eqref{mainlocal} only in the case that 
$g$ is constant at scale $\delta^{1/2}$, in which case there is no ambiguity.} 
\begin{equation}\label{mainlocal}
\| \widehat{g {\rm d} \sigma} \|_{L^{\frac{2n}{n-1}}(B(0,\delta^{-1}))} 
\leq C_n \left\| \left(\sum_\alpha | \widehat{g_\alpha {\rm d} \sigma} |^2\right)^{1/2} \right\|_{L^{\frac{2n}{n-1}}(B(0,\delta^{-1}))} 
\end{equation}
where  $g$ is a smooth function defined on $\mathbb{S}^{n-1}$, $\sigma$ is the Lebesgue measure on $\mathbb{S}^{n-1}$, 
$g_\alpha = g \chi_{E_\alpha}$, the $E_\alpha$ are spherical caps of radius $\delta^{1/2}$ decomposing $\mathbb{S}^{n-1}$ 
and where
$$\widehat{h {\rm d} \sigma}(x) = \int_{\mathbb{S}^{n-1}} h(\xi) e^{-2 \pi i x \cdot \xi} {\rm d} \sigma(\xi)$$
is the Fourier transform of the density $h {\rm d} \sigma$. Note that the inverse Fourier transform 
of $\widehat{g {\rm d} \sigma} \, \Phi( \delta \cdot)$ is essentially supported in a $\delta$ neighbourhood 
of $\mathbb{S}^{n-1}$, and that we may assume $g$ is constant at scale $\delta$.

\medskip
\noindent
The appropriate maximal function in this setting is the Kakeya maximal function which is defined as follows:
let $N \gg1 $ be a large parameter, and for $f$ defined on $\mathbb{R}^n$ and $\omega \in \mathbb{S}^{n-1}$
let 
$$M_N f(\omega) = \sup_{T \parallel \omega} \frac{1}{|T|} \int_T |f|$$
where the $\sup$ is taken over all tubes $T$ of dimensions $1 \times \dots \times 1 \times N$ whose axis 
is parallel to $\omega$. It is conjectured that 
\begin{equation}\label{maximal}
\|M_N f \|_{L^n(\mathbb{S}^{n-1})} \leq C_n \left( \log N \right)^\frac{n-1}{n} \|f \|_{L^n(\mathbb{R}^{n})}
\end{equation}
and this is also known to be true in two dimensions by work of Fefferman \cite{F} but is open in higher 
dimensions. Once again, \eqref{mainlocal} and \eqref{maximal} fit into a standard machine which can be used to 
establish, {\em inter alia}, restriction estimates such as
\begin{equation}\label{epr}
\| \widehat{g {\rm d} \sigma} \|_{L^{\frac{2n}{n-1}}(B(0,R))} \leq 
C_n  \left( \log R \right)^\frac{n-1}{2n}\|g\|_{L^{\frac{2n}{n-1}}(\mathbb{S}^{n-1})}.
\end{equation}

\medskip
\begin{proposition}\label{princ}
Suppose that \eqref{mainlocal} holds\footnote{for functions $g$ constant at scale $\delta^{1/2}$}, that is
$$\| \widehat{g {\rm d} \sigma} \|_{L^{\frac{2n}{n-1}}(B(0,\delta^{-1}))} 
\leq C_n \left\| \left(\sum_\alpha | \widehat{g_\alpha {\rm d} \sigma} |^2\right)^{1/2} \right\|_{L^{\frac{2n}{n-1}}(B(0,\delta^{-1}))}.$$
Then \eqref{maximal} holds, that is
$$\|M_N f \|_{L^n(\mathbb{S}^{n-1})} \leq C_n \left( \log N \right)^\frac{n-1}{n} \|f \|_{L^n(\mathbb{R}^{n})},$$ 
and hence 
$$\| \widehat{g {\rm d} \sigma} \|_{L^{\frac{2n}{n-1}}(B(0,R))} \leq 
C_n  \left( \log R \right)^\frac{n-1}{2n}\|g\|_{L^{\frac{2n}{n-1}}(\mathbb{S}^{n-1})}$$
also holds.
\end{proposition}

\medskip
\noindent
The proof uses a two-scale technique and is similar to Fefferman's argument for the disc multiplier \cite{FD} and 
to the argument in \cite{BCSS} in which it was first noted that restriction estimates imply estimates for maximal functions,
and which was later treated more formally by Bourgain \cite{B}. Also
see the remark at the end for another argument by-passing the maximal
function.

\medskip
\noindent
Since we know that the Kakeya maximal conjecture implies the Nikodym maximal conjecture (see Theorem 4.10 of 
\cite{T}, which is in turn based upon \cite{CC}) we conclude that \eqref{mainlocal} and \eqref{RLP} (together) 
also imply the Bochner--Riesz conjecture \eqref{BR}.

\medskip
\noindent
{\bf Remark.} This work was done in the early 1990's but was never presented for publication. The author has communicated its essence
to a number of harmonic analysts over the years and its existence has apparently become known, resulting in occasional 
requests for it. It is hoped that this informal presentation will satisfy such demand.

\section{Sketch of proof of Proposition \ref{princ}}

\noindent
We shall establish \eqref{maximal} in a certain dual form, see for example \cite{Co} and \cite{Carb}. For 
$N$ fixed there are essentially $N^{n-1}$ `distinct' directions that a $1 \times \dots \times 1 \times N$
tube can occupy. Suppose we have a collection $\{T_\alpha\}$ of such tubes, one in each direction. Then the best constant
$A_r$ in the inequality 
\begin{equation}\label{cov}
\| \sum_\alpha c_\alpha\chi_{T_\alpha} \|_r \leq A_r N^{1/r} \left( \sum_\alpha |c_\alpha|^r \right)^{1/r}
\end{equation}
(over all such families of tubes) is equivalent to the $L^{r'} \to L^{r'}$ operator norm of $M_N$. 
The equivalent scaled version of \eqref{cov} for families of 
$\lambda \times \dots \times \lambda \times \lambda N$ 
tubes with distinct directions is
$$\| \sum_\alpha c_\alpha\chi_{T_\alpha} \|_r \leq A_r N^{1/r} \lambda^{n/r} \left( \sum_\alpha |c_\alpha|^r \right)^{1/r}.$$
It is easy to verify that if the tubes $\{T_\alpha\}$ form a bush (i.e. all pass through 
a common centre) then \eqref{cov} holds for $r = n/(n-1)$ in the optimal form
\begin{equation}\label{covspec}
\| \sum_\alpha c_\alpha \chi_{T_\alpha} \|_{\frac{n}{n-1}} 
\leq C_n \left( \log N \right)^{\frac{n-1}{n}} N^{\frac{n-1}{n}} \left( \sum_\alpha |c_\alpha|^{\frac{n}{n-1}} \right)^{\frac{n-1}{n}}.
\end{equation}

\begin{proof}
{\bf \em Warning:} The treatment is somewhat informal. 

\medskip
\noindent
Let $T_\alpha$ be a collection of $\delta^{-1/2} \times \dots \times \delta ^{-1/2} \times \delta^{-1}$ tubes in 
$\mathbb{R}^n$ (so $\lambda = N = \delta^{-1/2}$ here), one in each of the $\delta^{-(n-1)/2}$ directions. By the 
remarks above, it suffices to show the appropriately scaled form of \eqref{cov} which is in this case
\begin{equation}\label{covscaled}
\| \sum_\alpha c_\alpha \chi_{T_\alpha} \|_{\frac{n}{n-1}} 
\leq C_n \left( \log N \right)^{\frac{n-1}{n}} N^{\frac{(n+1)(n-1)}{n}} \left( \sum_\alpha |c_\alpha|^{\frac{n}{n-1}} \right)^{\frac{n-1}{n}}.
\end{equation}

\noindent
We may assume that all the $T_\alpha$ are supported in a ball $B(0, \delta^{-1})$ in $\mathbb{R}^n$ and that $T_\alpha$ is the 
translate by $\lambda_\alpha$ of the tube passing through the origin in the $\alpha$'th direction. Note that $|\lambda_\alpha| 
\lesssim \delta^{-1}$. We may assume that $c_\alpha \geq 0$.

\medskip
\noindent
Let $\{E_\alpha\}$ be $\delta^{1/2}$-cells on $\mathbb{S}^{n-1}$ and let
$$ g(\xi) = \sum_{\alpha} \pm c_\alpha^{1/2} e^{i \lambda_\alpha \cdot \xi} \phi_\alpha(\xi)$$ 
where $\phi_\alpha$ is a smooth bump function associated to the cell $E_\alpha$.
Standard stationary phase calculations give
\begin{equation}\label{FT}
\left| \left( e^{i \lambda_\alpha \cdot \xi} \phi_\alpha (\xi){\rm d} \sigma (\xi)\right)^\wedge (x)\right| 
\gtrsim \delta^{(n-1)/2} \chi_{T_\alpha}(x).
\end{equation}
Since  $|\lambda_\alpha| \lesssim \delta^{-1}$ we have that $g$ is roughly constant at scale $\delta$ and so we 
can further decompose $g$ as 
$$ g = \sum_\beta a_\beta c_{\alpha(\beta)}^{1/2} \psi_\beta$$
where $\psi_{\beta}$ are smooth bump functions associated to a decomposition of the sphere into $\delta$-cells $F_\beta$, 
where $\alpha(\beta)$ is the $\alpha$ such that $F_\beta \subseteq E_\alpha$ and  where $|a_\beta| \sim 1$.

\medskip
\noindent
Now we apply our assumption \eqref{mainlocal} -- with $\delta^2$ now playing the role of 
$\delta$ in \eqref{mainlocal} -- to obtain
$$\| \widehat{g {\rm d} \sigma} \|_{L^{\frac{2n}{n-1}}(B(0,\delta^{-2}))} 
\leq C_n \left\| \left(\sum_\beta c_{\alpha(\beta)}| \widehat{\psi_\beta {\rm d} \sigma} |^2\right)^{1/2} \right\|_{L^{\frac{2n}{n-1}}(B(0,\delta^{-2}))}. $$
The main contribution to $| \widehat{\psi_\beta {\rm d} \sigma} |$ is given by $ \delta^{(n-1)}\chi_{R_\beta}$
where $R_\beta$ is a $\delta^{-1} \times \dots \delta^{-1} \times \delta^{-2}$ tube {\em passing through the origin}.
Hence, ignoring lower order contributions,
$$\left\|  \sum_{\alpha} \pm c_\alpha^{1/2} \left(e^{i \lambda_\alpha \cdot \xi} \phi_\alpha(\xi){\rm d} \sigma (\xi)\right)^\wedge \right\|_{L^{\frac{2n}{n-1}}(B(0,\delta^{-2}))} $$
$$\leq 
C_n  \delta^{(n-1)} \left\| \left(\sum_\beta c_{\alpha(\beta)} \chi_{R_\beta}\right)^{1/2} \right\|_{L^{\frac{2n}{n-1}}(B(0,\delta^{-2}))}. $$

\medskip
\noindent
Next we use Khintchine's inequality and \eqref{FT} to deduce that
$$ \delta^{(n-1)/2} \left\| \left(\sum_\alpha c_{\alpha} \chi_{T_\alpha}\right)^{1/2} \right\|_{L^{\frac{2n}{n-1}}(B(0,\delta^{-2}))} 
\leq C_n\delta^{(n-1)} \big\| \left(\sum_\beta c_{\alpha(\beta)} \chi_{R_\beta}\right)^{1/2} \big\|_{L^{\frac{2n}{n-1}}(B(0,\delta^{-2}))} $$
or equivalently
$$ \left\| \sum_\alpha c_{\alpha} \chi_{T_\alpha} \right\|_{\frac{n}{n-1}} 
\leq C_n\delta^{(n-1)} \left\| \sum_\beta c_{\alpha(\beta)} \chi_{R_\beta} \right\|_{\frac{n}{n-1}}.$$

\medskip
\noindent
Now the tubes $R_\beta$ all pass through the origin, so by the scaled version of the remark immediately preceeding the proof we 
have
$$\left\| \sum_\beta d_{\beta} \chi_{R_\beta} \right\|_{\frac{n}{n-1}} \leq 
C_n \left( \log \frac{1}{\delta}\right)^{\frac{n-1}{n}} \delta^{\frac{-(n+1)(n-1)}{n}} \left(\sum_\beta |d_\beta|^{\frac{n}{n-1}}\right)^{\frac{n-1}n}.$$
Taking $d_\beta = c_{\alpha(\beta)}$ and noting that 
$$ \sum_\beta |d_\beta|^{\frac{n}{n-1}} = \delta^{-(n-1)/2} \sum_\alpha |c_\alpha|^{\frac{n}{n-1}}$$
we obtain

$$ \left\| \sum_\alpha c_{\alpha} \chi_{T_\alpha} \right\|_{\frac{n}{n-1}} 
\leq C_n \left( \log \frac{1}{\delta} \right)^{\frac{n-1}{n}} \delta^{(n-1)} \delta^{\frac{-(n-1)^2}{2n}} \delta^{\frac{-(n+1)(n-1)}{n}} 
\left(\sum_\alpha |c_\alpha|^{\frac{n}{n-1}}\right)^{\frac{n-1}{n}}$$

$$ =  C_n \left( \log \frac{1}{\delta} \right)^{\frac{n-1}{n}} \delta^{\frac{-(n+1)(n-1)}{2n}} 
\left(\sum_\alpha |c_\alpha|^{\frac{n}{n-1}}\right)^{\frac{n-1}{n}}.$$
\end{proof}

\medskip
\noindent
{\bf Remark.}
Note that we did not need the full force of hypothesis \eqref{mainlocal} for general $g$, only for those 
$g$ constant on $\delta^{1/2}$-cells. 

\medskip
\noindent
Indeed, suppose we assume \eqref{mainlocal} for such $g$ constant on $\delta^{1/2}$ cells $E_\alpha$. We can then 
conclude inequality
\eqref{epr} directly without passing through the maximal function. To see this, we begin by observing that 
in order to prove inequality \eqref{epr} we may assume that $g$ is a step function on $\mathbb{S}^{n-1}$
which is constant at scale $R^{-1}$ (this is the effect of restricting attention to $B(0,R)$). So, since 
$B(0,R) \subseteq B(0,R^2)$, 
in order to prove \eqref{epr} for such $g$ it suffices to prove the ostensibly stronger\footnote{but not really, as 
$ \widehat{g {\rm d} \sigma}$ is pretty much supported in $B(0,R)$} inequality
$$\| \widehat{g {\rm d} \sigma} \|_{L^{\frac{2n}{n-1}}(B(0,R^2))} \leq 
C_n  \left( \log R \right)^\frac{n-1}{2n}\|g\|_{L^{\frac{2n}{n-1}}(\mathbb{S}^{n-1})}$$
for such $g$. But if we relabel $R^2$ as $R$, we see it therefore suffices to prove that inequality \eqref{epr}
holds for $g$ constant at scale $R^{-1/2}$; that is, we may assume $g = \sum_\alpha c_\alpha \phi_\alpha$. 
For such $g$ we are assuming that we have inequality \eqref{mainlocal} with $\delta = 1/R$ 
and so 
$$ \| \widehat{g {\rm d} \sigma} \|_{L^{\frac{2n}{n-1}}(B(0,R))} \leq  
C_n \left\| \left(\sum_\alpha | c_\alpha|^2 | \widehat{\phi_\alpha {\rm d} \sigma} |^2\right)^{1/2} \right\|_{L^{\frac{2n}{n-1}}(B(0,R))}$$ 
where the main contributions of the terms  $|\widehat{\phi_\alpha {\rm d} \sigma}|$ come from 
$R^{1/2} \times \dots \times R^{1/2} \times R$-tubes centred at the origin, and which we may therefore calculate and estimate 
directly. Indeed, for such $g$ constant at scale $R^{-1/2}$, we have\footnote{There are similar expressions valid for other 
values of $q$ and $r$.}
$$ \left\| \left(\sum_\alpha | c_\alpha|^2 | \widehat{\phi_\alpha {\rm d} \sigma} |^2\right)^{1/2} \right\|_{L^{\frac{2n}{n-1}}(B(0,R))}
\approx \left(\int_{{R^{-1/2}}}^1 \int_{\mathbb{S}^{n-1}} \left(\mathbb{A}_t |g|^2\right)^{n/(n-1)} {\rm d} \sigma \frac{{\rm d}t}{t}\right)^{(n-1)/2n}$$
where $\mathbb{A}_t$ denotes a local average on the sphere at scale $t$. Now by the triangle inequality we have
$$ \left(\int_{\mathbb{S}^{n-1}} \left(\mathbb{A}_t |g|^2\right)^{n/(n-1)} {\rm d} \sigma \right)^{(n-1)/n} \leq \||g|^2\|_{n/n-1}$$ 
and so 
$$ \left(\int_{R^{-1/2}}^{1} \int_{\mathbb{S}^{n-1}} \left(\mathbb{A}_t |g|^2(\omega)\right)^{n/(n-1)} {\rm d} \sigma(\omega) \frac{{\rm d}t}{t}\right)^{(n-1)/2n}
\leq C_n  \left( \log R \right)^{\frac{n-1}{2n}} \|g\|_{L^{\frac{2n}{n-1}}(\mathbb{S}^{n-1})}.$$


\medskip
\noindent
See also the articles \cite{CoM}, \cite{MVV1} and \cite{MVV2} for related discussions.

\bibliographystyle{plain}

\end{document}